\documentclass[11pt,reqno]{amsart}
\usepackage{amsmath,amsfonts, amssymb, amsthm}
  \usepackage{paralist}
  \usepackage{graphics} 
  \usepackage{epsfig} 
\usepackage{graphicx}  \usepackage{epstopdf}
 \usepackage[colorlinks=true]{hyperref}
\hypersetup{urlcolor=blue, citecolor=red}

\newtheorem{theorem}{Theorem}[section]
\newtheorem{corollary}{Corollary}[section]

\newtheorem{lemma}[theorem]{Lemma}
\newtheorem{proposition}{Proposition}[section]
\newtheorem{conjecture}{Conjecture}[section]

\theoremstyle{definition}

\newtheorem{remark}{Remark}[section]

\def\pmod #1{\ ({\rm{mod}}\ #1)}
\def\Z{\Bbb Z}
\def\N{\Bbb N}

\def\Q{\Bbb Q}

\def\R{\Bbb R}
\def\C{\Bbb C}

\def\l{\left}
\def\r{\right}
\def\bg{\bigg}
\def\({\bg(}
\def\){\bg)}
\def\t{\text}
\def\f{\frac}
\def\mo{{\rm{mod}\ }}
\def\pmod#1{\ (\mo\ #1)}

\def\ls{\leq}
\def\gs{\geq}

\def\bi{\binom}
\def\al{\alpha}

\def\eq{\equiv}

\def\Proof{\noindent{\it Proof}}

\begin{document}
\hbox{Nanjing Univ. J. Math. Biquarterly 41 (2024), no. 1, 1--33.}
\medskip

\title[Series of the type $\sum_{k=0}^\infty(ak+b)x^k/\binom{mk}{nk}$]
      {Evaluations of some series of the type $\sum_{k=0}^\infty(ak+b)x^k/\binom{mk}{nk}$}
\author[Zhi-Wei Sun]{Zhi-Wei Sun}


\address{Department of Mathematics, Nanjing
University, Nanjing 210093, People's Republic of China}
\email{{\tt zwsun@nju.edu.cn}
\newline\indent
{\it Homepage}: {\tt http://maths.nju.edu.cn/\lower0.5ex\hbox{\~{}}zwsun}}

\keywords{Series for $\pi$, logarithmic function, central binomial coefficients, combinatorial identities, beta function.
\newline \indent The initial version of this paper was posted to arXiv in 2022 with the ID 2204.08275.
\newline \indent Supported by the National Natural Science Foundation of China (grant no. 12371004).}

\begin{abstract}
In this paper, via the beta function we evaluate some series of the type
$\sum_{k=0}^\infty(ak+b)x^k/\binom{mk}{nk}$.
For example, we prove that
$$\sum_{k=0}^\infty\frac{(49k+1)8^k}{3^k\binom{3k}k}=81+16\sqrt3\,\pi
\ \ \t{and}\ \
\sum_{k=0}^\infty\frac{10k-1}{\binom{4k}{2k}}=\frac{4\sqrt 3}{27}\pi.$$
 We also establish the following efficient formula
for computing $\log n$ with $1<n\ls 85/4$:
\begin{align*}
&\sum_{k=0}^\infty\frac{(2(n^2+6n+1)^2(n^2-10n+1)k+P(n))(n-1)^{4k}}
{(-n)^k(n+1)^{2k}\binom{4k}{2k}}\\ \ \ &=6n(n+1)(n-1)^3\log n-32n(n+1)^2(n^2-4n+1),
\end{align*}
where
$$P(n):=n^6-58n^5+159n^4+52n^3+159n^2-58n+1.$$
In addition, we pose some conjectures on series whose summands involve $\bi{2k}k/(\bi{3k}k\bi{6k}{3k})\ (k\in\N)$.
\end{abstract}
\maketitle

\section{Introduction}

For any real number $x$ with $|x|<2$, {\tt Mathematica} yields
\begin{equation}\label{+arcsin}\f14\sum_{k=0}^\infty\f{x^{2k}}{\bi {2k}k}=\f{\sqrt{4-x^2}+x\arcsin (x/2)}{(4-x^2)\sqrt{4-x^2}}\end{equation}
and
\begin{equation}\label{-arcsin}
\f14\sum_{k=0}^\infty\f{(-1)^kx^{2k}}{\bi{2k}k}=\f{1}{4+x^2}-\f{x\,\mathrm{arcsinh}(x/2)}{(4+x^2)\sqrt{4+x^2}},
\end{equation}
where
$$\mathrm{arcsinh}\ t=\sum_{n=0}^\infty\f{\bi{2n}nt^{2n+1}}{(2n+1)(-4)^n}=\log(t+\sqrt{t^2+1})$$
is the inverse hyperbolic sine function.
It is also known that
$$\sum_{k=1}^\infty\f{x^{2k}}{k^2\bi{2k}k}=2\arcsin^2\f x2$$
for any $x\in\R$ with $|x|<2$ (see, e.g., \cite{BC}).

Using \eqref{+arcsin} and \eqref{-arcsin} and their derivatives, we can easily deduce that
\begin{equation}\sum_{k=0}^\infty\f{2(2n+1)^2k+3}{(-n(n+1))^k\bi{2k}k}=-\f{2n(n+1)}{2n+1}\log\l(1+\f1n\r)
\end{equation}
if $n<-(1+\sqrt2)/2$ or $n>(\sqrt2-1)/2$.
When $n=1/4$ this yields the identity
$$\sum_{k=0}^\infty\f{(3k+2)16^k}{(-5)^k\bi{2k}k}=-\f5{18}\log 5.$$
In contrast, if $n>1$ or $n<-1$ then
$$\log\l(1+\f1n\r)=\sum_{k=1}^\infty\f{(-1)^{k-1}}{kn^k}$$
by the Taylor series.

Series for $\pi$ are particularly interesting. For Ramanujan-type series, one may consult S. Cooper \cite[Chapter 14]{Co}. For few double series for $\pi$, one may consult C. Wei \cite{Wei}.
In 1974 R. W. Gosper announced the new identity
\begin{equation}\label{Gosper}\sum_{k=0}^\infty\f{25k-3}{2^k\bi{3k}k}=\f{\pi}2,
\end{equation}
which was later used by F. Bellard \cite{Be} to find an algorithm for computing the $n$th decimal of $\pi$ without calculating the earlier ones. Inspired by Gosper's identity, in 2003 Bellard \cite{Be} discovered the identity
$$\pi=\f1{740025}\(\sum_{k=1}^\infty\f{3P(k)}{2^{k-1}\bi{7k}{2k}}-20379280\),$$
where
\begin{align*}P(n)&=-885673181k^5+3125347237k^4-2942969225k^3
\\&\ \ +1031962795k^2-196882274k+10996648;
\end{align*}
he used this identity to set his world record of computing the $10^{11}$ binary digit of $\pi$. Moreover,
G. Almkvist, C. Krattenthaler
and J. Petersson \cite{AKP} gave a proof of Gosper's identity and
found 12 new identities of the type
$$\pi=\sum_{k=0}^\infty\f{P(k)}{a^k\bi{mk}{nk}},$$
where $P(x)\in\Q[x]$, and $(m,n,a,\deg P)$ is among the ordered quadruples
\begin{gather*}(8,4,-4,4),\,(10,4,4,8),\,(12,4,-4,8),\,(16,8,16,8),
\\(24,12,-64,12),\,(32,16,256,16),(40,20,-2^{10},20),\,(48,24,2^{12},24),
\\(56,28,-2^{14},28),(64,32,2^{16},32),\,(72,36,-2^{18},36),\,(80,40,2^{20},40).
\end{gather*}
For example, \cite[Example 2]{AKP} gives the identity
$$\pi=\f1{105^2}\sum_{k=0}^\infty\f{P(k)}{\bi{8k}{4k}(-4)^k},$$
where
$$P(k)=-89286+3875948k-34970134k^2+110202472k^3-115193600k^4.$$

By Stirling's formula,
$$k!\sim\sqrt{2\pi k}\l(\f ke\r)^k\ \t{as}\ k\to+\infty.$$
If $m>n>0$ are integers, then
$$\bi{mk}{nk}=\f{(mk)!}{(nk)!((m-n)k)!}\sim\f{\sqrt m}{\sqrt{2\pi n(m-n)k}}\l(\f {m^m}{n^n(m-n)^{m-n}}\r)^k$$
as $k\to+\infty.$

 In this paper we evaluate some series of the type
 $$\sum_{k=0}^\infty(ak+b)\f{x^k}{\bi{mk}{nk}},$$
 where $m>n>0$ are integers and $a,b,x$ are real numbers with
 $$|x|<\f{m^m}{n^n(m-n)^{m-n}}.$$

 Note that $$\bi{3k}k\sim\f{\sqrt3}{2\sqrt{k\pi}}\l(\f{27}4\r)^k\ (k\to+\infty).$$
Thus, for any real number $x_0$ with $-27/4<x_0<27/4$ the series
$$\sum_{k=0}^\infty\f{x_0^k}{\bi{3k}k}\ \ \t{and}\ \ \sum_{k=0}^\infty\f{kx_0^k}{\bi{3k}k}$$
converge absolutely.
Let
\begin{equation}\label{c}c=\f32\l((1+\sqrt2)^{1/3}-(1+\sqrt2)^{-1/3}\r)=0.8941\cdots.
\end{equation} For $f(x)=x^3-x_0(x-1)$,
clearly $f(-3)=-27+4x_0<0$ and $$f(c)=c^3+(1-c)x_0=\f{27}4(1-c)+(1-c)x_0=(1-c)\l(x_0+\f{27}4\r)>0.$$
So there is a real number $-3<x<c$ such that $f(x)=0$ and hence $x_0=x^3/(x-1)$.
Moreover, such $x\in(-3,c)$ can be found by solving the cubic equation $x^3=x_0(x-1)$.

 Now we state our first theorem.

 \begin{theorem}\label{C(3k,k)} Let $-3<x<c$ with $c$ given by \eqref{c}.

 {\rm (i)} We have
 \begin{equation}\label{arctan}
 \begin{aligned}&\sum_{k=1}^\infty\f{((2x-3)^2k+2x^2+2x-3)x^{3k}}{(x-1)^k\bi{3k}k}
 \\=\ &-2x^3\f{x+7}{(x+3)^2}
 +\f{8x^2(x-1)q(x)}{(x+3)^2\sqrt{(1-x)(3+x)}},
 \end{aligned}
 \end{equation}
 where
 $$q(x)=\begin{cases}\arctan\f x{x+2}\sqrt{\f{3+x}{1-x}}&\t{if}\ -2<x<1,
 \\-\f{\pi}2&\t{if}\ x=-2,
 \\\arctan\f x{x+2}\sqrt{\f{3+x}{1-x}}-\pi&\t{if}\ -3<x<-2.
 \end{cases}$$

 {\rm (ii)} We have
 \begin{equation}\label{log}
 \sum_{k=0}^\infty\f{(s(x)k+t(x))x^{3k}}{(x-1)^k\bi{3k}k}
 =12x^2(1-x)\log(1-x)-27(1-x)(x^2-6x+3),
 \end{equation}
 where
 \begin{equation}\label{s(x)}
 s(x)=(x+3)(2x-3)^2(x^2-12x+9)=4x^5-48x^4+9x^3+351x^2-567x+243
 \end{equation}
 and
 \begin{equation}\label{t(x)}
 t(x)=2x^5-48x^4+69x^3-189x^2+243x-81.
 \end{equation}
 \end{theorem}

 Since
 $$\left|\begin{matrix}
(2x-3)^2&2x^2+2x-3\\s(x)&t(x)\end{matrix}\right|  =-4x(2x-3)^5\not=0$$
 for any $x\in(-3,c)$ with $x\not=0$, a suitable combination of the two parts of Theorem \ref{C(3k,k)} yields the values of
 $$\sum_{k=0}^\infty\f{x^{3k}}{(x-1)^k\bi{3k}k}\ \ \t{and}\ \
 \sum_{k=0}^\infty\f{kx^{3k}}{(x-1)^k\bi{3k}k}$$
 for all $x\in(-3,c)$.

 \begin{corollary} Whenever  $-3<x<c$, we have
 \begin{equation}\label{explicit}
\begin{aligned}
\sum_{k=0}^\infty\f{x^{3k}}{(x-1)^k\bi{3k}k}=&\ \f{27(1-x)}{(x+3)(2x-3)^2}+\f{3x(x-1)}{(2x-3)^3}\log(1-x)
\\&\ +\f{2x(x-1)(x^2-12x+9)q(x)}{(x+3)(2x-3)^3\sqrt{(1-x)(x+3)}}
\end{aligned}\end{equation}
\end{corollary}
\begin{remark} Clearly, $s(x)\times$\eqref{arctan} minus $(2x-3)^2\times$\eqref{log} yields the equality
\eqref{explicit}
Note that N. Batir \cite[(3.3)]{Ba} gave a very complicated formula for
$$\sum_{k=1}^\infty\f{x_0^k}{k^r\bi{3k}k}\ \ \l(-\f{27}4<x_0<\f{27}4\r).$$
 \end{remark}

 Observe that \eqref{arctan} in the case $x=-1$ gives Gosper's identity \eqref{Gosper}.

 Putting $x=-2$ in \eqref{arctan} and $x=1/n$ in \eqref{log}, we obtain the following corollary.

 \begin{corollary} {\rm (i)} We have
 \begin{equation}\sum_{k=0}^\infty\f{(49k+1)8^k}{3^k\bi{3k}k}=81+16\sqrt3\,\pi.
 \end{equation}

 {\rm (ii)} If $n<-1/3$ or $n>1/c=1.11843\cdots$, then
 \begin{equation}\sum_{k=0}^\infty\f{a_nk-b_n}{((1-n)n^2)^k\bi{3k}k}=3n^2(n-1)\l(4\log\l(1-\f1n\r)
 -9(3n^2-6n+1)\r),
 \end{equation}
 where
 $$a_n=(3n+1)(3n-2)^2(9n^2-12n+1)$$
 and
 $$b_n=81n^5-243n^4+189n^3-69n^2+48n-2.$$
 In particular,
 \begin{align}\sum_{k=0}^\infty\f{275k-158}{2^k\bi{3k}k}&=6\log2-135,
 \\\sum_{k=0}^\infty\f{728k-17}{(-4)^k\bi{3k}k}&=-54-24\log2,
 \\\sum_{k=0}^\infty\f{(1813k-2707)8^k}{3^k\bi{3k}k}&=9(16\log3-171),
 \\\sum_{k=0}^\infty\f{5635k-1156}{(-18)^k\bi{3k}k}&=54\log\f23-1215,
 \\\sum_{k=0}^\infty\f{63050k-15959}{(-48)^k\bi{3k}k}&=72\l(4\log\f 34-225\r),
 \\\sum_{k=0}^\infty\f{112216k-30847}{(-100)^k\bi{3k}k}&=300\log\f 45-31050,
 \\\sum_{k=0}^\infty\f{615296k-176777}{(-180)^k\bi{3k}k}&=270\l(4\log\f 56-657\r),
 \\\sum_{k=0}^\infty\f{710809k-209926}{(-294)^k\bi{3k}k}&=441\l(2\log\f 67-477\r),
 \end{align}
 \begin{align}
 \sum_{k=0}^\infty\f{2910050k-875807}{(-448)^k\bi{3k}k}&=672\l(4\log\f78-1305\r),
 \\\sum_{k=0}^\infty\f{2721250k-830317}{(-648)^k\bi{3k}k}&=972\l(2\log\f89-855\r),
 \\\sum_{k=0}^\infty\f{9490712k-2926289}{(-900)^k\bi{3k}k}&=1350\l(4\log\f 9{10}-2169\r),
 \\\sum_{k=0}^\infty\f{7825423k-2432776}{(-1210)^k\bi{3k}k}&=1815\l(2\log\f {10}{11}-1341\r).
 \end{align}
 \end{corollary}

  Let $x\in(-16,16)$. By induction, we have
 \begin{align*}&\sum_{k=1}^n\f{x^{k}}{\bi{4k}{2k}}\l(-\f 6k+x+32+2k(x-16)\r)
=-x+\f{(2n+1)x^{n+1}}{\bi{4n}{2n}}.
\end{align*}
which tends to $-x$ as $n\to+\infty$. Thus, if we know the values of
$$\sum_{k=0}^\infty\f{x^k}{\bi{4k}{2k}}\ \t{and}\ \sum_{k=0}^\infty\f{kx^k}{\bi{4k}{2k}}$$
then the value of $\sum_{k=1}^\infty x^{k}/(k\bi{4k}{2k})$ is also determined.

 For any $x\in\R$ with $x>1$ or $x<0$, we define
 \begin{equation}\label{R(x)}R(x):=\sqrt x\,\mathrm{arctanh}\f1{\sqrt x}=\begin{cases}\f{\sqrt x}2\log\f{\sqrt x+1}{\sqrt x-1}
 &\t{if}\ x>1,
 \\ \sqrt{|x|}\arctan\f1{\sqrt{|x|}}&\t{if}\ x<0,\end{cases}
 \end{equation}
 where $\mathrm{arctanh}\,t$ is the inverse hyperbolic tangent function. Note that for $x>1$ or $x<-1$,
 we have
 $$R(x)=\sqrt x\sum_{k=0}^\infty\f{(1/\sqrt x)^{2k+1}}{2k+1}=\sum_{k=0}^\infty\f {x^{-k}}{2k+1}.$$

Now we are ready to state our second theorem which involves the binomial coefficients $\bi{4k}{2k}$ with $k\in\N$.

 \begin{theorem}\label{Th1.1} For any $x>1/4$, we have
 \begin{equation}\label{comput-arccot}\sum_{k=0}^\infty\f{2(4x+1)k-2x+1}{x^{2k}\bi{4k}{2k}}
 =\f{8x^2}{(4x-1)^2}\l(\f3{\sqrt{4x-1}}\mathrm{arccot}\sqrt{4x-1}-4x+4\r)
 \end{equation}
 and
 \begin{equation}\label{comput-arccoth}\sum_{k=0}^\infty\f{2(4x-1)k-2x-1}{x^{2k}\bi{4k}{2k}}
 =\f{8x^2}{(4x+1)^2}\l(\f{3R(4x+1)}{4x+1}-4x-4\r).
 \end{equation}
 Consequently,
 \begin{align}\label{new-1}\sum_{k=0}^\infty\f{10k-1}{\bi{4k}{2k}}&=\f{4\sqrt 3}{27}\pi,
 \\\label{New1}\sum_{k=0}^\infty\f{k4^k}{\bi{4k}{2k}}&=\f{3\pi+8}{12},
 \\\label{new-9}\sum_{k=0}^\infty\f{(14k+1)9^k}{\bi{4k}{2k}}&=24\pi\sqrt3+64,
 \\\label{new-9/4}\sum_{k=0}^\infty\f{(22k-1)9^k}{4^k\bi{4k}{2k}}&=\f{32}{25}\l(4+\f{27}{\sqrt{15}}\arctan\sqrt{\f35}\r),
 \\\label{New3}\sum_{k=0}^\infty\f{14k-5}{4^k\bi{4k}{2k}}&=\f{16}{81}(\log2-24).
 \end{align}
 \end{theorem}
 \begin{remark} (a) As $$\left|\begin{matrix}
2(4x+1)&-2x+1\\2(4x-1)&-2x-1\end{matrix}\right|  =-24x\not=0$$
 for all $x>1/4$, combining \eqref{comput-arccot} and \eqref{comput-arccoth} we have actually determined the values of
 $$\sum_{k=0}^\infty\f{x_0^{k}}{\bi{4k}{2k}}\ \ \t{and}\ \
 \sum_{k=0}^\infty\f{kx_0^{k}}{\bi{4k}{2k}}$$
 for all $x_0\in(0,16)$.

 (b) In contrast with \eqref{New1}, for any prime $p>3$ we conjecture the congruences
 \begin{equation*}p\sum_{k=0}^{(p-1)/2}\f{k4^k}{\bi{4k}{2k}}\eq\f12\l(\f{-1}p\r)-\f p6\pmod{p^2}.
 \end{equation*}
 and
 \begin{equation*}p\sum_{k=0}^{p-1}\f{k4^k}{\bi{4k}{2k}}\eq\f23\l(\f{-1}p\r)+\f29p\pmod{p^2},
 \end{equation*}
 where $(\f{\cdot}p)$ demotes the Legendre symbol.
 \end{remark}

For $x_0\in(0,16)$, how to evaluate
$$\sum_{k=0}^\infty\f{(-x_0)^k}{\bi{4k}{2k}}\ \ \t{and}\ \ \sum_{k=0}^\infty\f{k(-x_0)^k}{\bi{4k}{2k}}?$$
If we take
$$x=\f12+\sqrt{\f 4{x_0}+\f14}>\f12+\sqrt{\f14+\f14}=\f{1+\sqrt2}2,$$
then
$$-x_0=\f4{x(1-x)}.$$

 Now we state our third theorem.

 \begin{theorem} \label{C(4k,2k)} If $x>(1+\sqrt2)/2$ or $x<(1-\sqrt2)/2$, then
 \begin{equation}\label{gen}
 \begin{aligned}&\sum_{k=0}^\infty\f{(2(2x-1)^2(2x-3)k-(4x^3-16x^2+7x+6))4^k}{(x(1-x))^k\bi{4k}{2k}}
 \\&\ \ =(1-x)\l(3R(x)+4x(x-3)\r)
 \end{aligned}
 \end{equation}
 and
 \begin{equation}\label{gen-dual}
 \begin{aligned}&\sum_{k=0}^\infty\f{(2(2x-1)^2(2x+1)k-(4x^3+4x^2-13x-1))4^k}{(x(1-x))^k\bi{4k}{2k}}
 \\&\ \ =-x\l(3R(1-x)+4(x-1)(x+2)\r).
 \end{aligned}
 \end{equation}
 \end{theorem}
 \begin{remark} As $$\left|\begin{matrix}
2(2x-1)^2(2x-3)&-(4x^3-16x^2+7x+6)\\2(2x-1)^2(2x+1)&-(4x^3+4x^2-13x-1)\end{matrix}\right|  =-6(2x-1)^5\not=0$$
  if $x>(1+\sqrt2)/2$ or $x<(1-\sqrt2)/2$, combining the two identities in Theorem \ref{C(4k,2k)} we have actually determined the values of
 $$\sum_{k=0}^\infty\f{(-x_0)^{k}}{\bi{4k}{2k}}\ \ \t{and}\ \
 \sum_{k=0}^\infty\f{k(-x_0)^{k}}{\bi{4k}{2k}}$$
 for all $x_0\in(0,16)$.
 \end{remark}

 \begin{corollary} We have
 \begin{equation}\label{New2}\sum_{k=0}^\infty\f{(30k-7)(-2)^k}{\bi{4k}{2k}}=-\f{3\pi+64}{6}.
 \end{equation}
 \end{corollary}
 \Proof.  As $R(-1)=\arctan1=\pi/4$, clearly \eqref{New2} follows from
  \eqref{gen} with $x=-1$. \qed

\begin{remark} For any prime $p>3$, we conjecture the congruences
 \begin{equation*}p\sum_{k=0}^{(p-1)/2}\f{(30k-7)(-2)^k}{\bi{4k}{2k}}\eq-2\l(\f{-1}p\r)\pmod{p}
 \end{equation*}
 and
 \begin{equation*}p\sum_{k=0}^{p-1}\f{(30k-7)(-2)^k}{\bi{4k}{2k}}\eq-\f53\l(\f{-1}p\r)-\f{128}9p\pmod{p^2}.
 \end{equation*}
 \end{remark}

\begin{corollary} For
$$1<n<\f{\sqrt{(1+\sqrt2)/2}+1}{\sqrt{(1+\sqrt2)/2}-1}=21.2666866\ldots,$$
we have the following formula for $\log n$:
\begin{equation}\label{log n}\begin{aligned}
&\sum_{k=0}^\infty\f{(2(n^2+6n+1)^2(n^2-10n+1)k+P(n))(n-1)^{4k}}
{(-n)^k(n+1)^{2k}\bi{4k}{2k}}\\ \ \ &=6n(n+1)(n-1)^3\log n-32n(n+1)^2(n^2-4n+1),
\end{aligned}\end{equation}
where
$$P(n):=n^6-58n^5+159n^4+52n^3+159n^2-58n+1.$$
In particular,
 \begin{align}\label{log(2)}\sum_{k=0}^\infty\f{2890k-563}{(-18)^k\bi{4k}{2k}}&=-12(\log2+48),
 \\\sum_{k=0}^\infty\f{245k-17}{(-3)^k\bi{4k}{2k}}&=-24-\f 92\log3,
 \\\sum_{k=0}^\infty\f{(77326k+8951)81^k}{(-100)^k\bi{4k}{2k}}&=40(80-81\log4),
 \\\sum_{k=0}^\infty\f{(196k+73)64^k}{(-45)^k\bi{4k}{2k}}&=15(3-\log5),
 \\
 \sum_{k=0}^\infty\f{(245134k+181679)625^k}{(-294)^k\bi{4k}{2k}}&=84(1456-375\log6),
 \\\sum_{k=0}^\infty\f{(2645k+3517)81^k}{(-28)^k\bi{4k}{2k}}&=7(352-81\log 7),
 \\\sum_{k=0}^\infty\f{(127890k+316933)2401^k}{(-648)^k\bi{4k}{2k}}&=144(1584-343\log8),
 \\\sum_{k=0}^\infty\f{(1156k+7031)1024^k}{(-225)^k\bi{4k}{2k}}&=45(115-24\log 9),
 \\\sum_{k=0}^\infty\f{(51842k-3142679)6561^k}{(-1210)^k\bi{4k}{2k}}&=220(2187\log10-10736),
 \\
 \sum_{k=0}^\infty\f{(2209k-13421)625^k}{(-99)^k\bi{4k}{2k}}&=\f{99}2(125\log11-624),
 \\\sum_{k=0}^\infty\f{(2354450k-8037191)14641^k}{(-2028)^k\bi{4k}{2k}}&=312(3993\log12-20176),
 \\
 \sum_{k=0}^\infty\f{(19220k-46979)5184^k}{(-637)^k\bi{4k}{2k}}&=91(81\log13-413),
 \end{align}
 \begin{align}
 \sum_{k=0}^\infty\f{(3000515k-5794357)28561^k}{(-3150)^k\bi{4k}{2k}}&=420(2197\log14-11280),
 \\
 \sum_{k=0}^\infty\f{(118579k-190573)2401^k}{(-240)^k\bi{4k}{2k}}&=30(1029\log15-5312),
 \\\sum_{k=0}^\infty\f{(24174146k-33367199)50625^k}{(-4624)^k\bi{4k}{2k}}
 &=544(10125\log16-52496),
 \\
 \sum_{k=0}^\infty\f{(48020k-58117)16384^k}{(-1377)^k\bi{4k}{2k}}&=459(64\log17-333),
 \\
 \sum_{k=0}^\infty\f{(54371810k-58537799)83521^k}{(-6498)^k\bi{4k}{2k}}
 &=684(14739\log18-76912),
 \\\sum_{k=0}^\infty\f{(608923k-589327)6561^k}{(-475)^k\bi{4k}{2k}}&=\f{95}2(2187\log19-11440),
 \\\sum_{k=0}^\infty\f{(36377094k-31893853)130321^k}{(-8820)^k\bi{4k}{2k}}&=840(6859\log20-35952),
 \\\label{log(21)}\sum_{k=0}^\infty\f{(584756k-467339)40000^k}{(-2541)^k\bi{4k}{2k}}
 &=231(375\log21-1969),
 \end{align}
 and
 \begin{equation}\begin{aligned}&\sum_{k=0}^\infty\f{(661704134402k-517115569199)43046721^k}{2693140^k\bi{4k}{2k}}
 \\\quad&=60520\l(1594323\log\f{85}4-8374544\r).
 \end{aligned}
 \end{equation}
 \end{corollary}
 \Proof. Putting $x=(n+1)^2/(n-1)^2$ in \eqref{gen}, we get \eqref{log n}.
 Taking $n=2,\ldots,21,85/4$ in \eqref{log n} we immediately obtain the remaining identities. \qed

 \begin{remark} Note that our identities \eqref{log(2)}-\eqref{log(21)} provide series for
 $$\log 2,\ \ldots,\ \log21$$ which converge rapidly.
 The identity \eqref{log n} with $n=5/3,7/5,9/7$ yields the following examples:
 \begin{align}
 \sum_{k=0}^\infty\f{27869k-6203}{(-60)^k\bi{4k}{2k}}&=-15\l(416+3\log\f 53\r),
 \\\sum_{k=0}^\infty\f{115943k-27691}{(-315)^k\bi{4k}{2k}}&=-\f{105}2\l(528+3\log\f 75\r),
 \\\sum_{k=0}^\infty\f{2016125k-491747}{(-1008)^k\bi{4k}{2k}}&=-126\l(3904+3\log\f 97\r).
 \end{align}
 \end{remark}

In the next section we shall give an auxiliary proposition whose proof involves  the beta function
 $$B(a,b):=\int_0^1 x^{a-1}(1-x)^{b-1}dx\ \ \ \ \t{for}\ a>0\ \t{and}\ b>0.$$
 Our proofs of Theorems \ref{C(3k,k)}-\ref{C(4k,2k)}
 will be given in Sections 3--5 respectively.
 In Section 6, we pose some conjectural series whose summands involve $\bi{2k}k/(\bi{3k}k\bi{6k}k)$.

\section{An auxiliary proposition}
 \setcounter{equation}{0}
 \setcounter{conjecture}{0}
 \setcounter{theorem}{0}
 \setcounter{proposition}{0}

 \begin{lemma} \label{Lem2.1} For any complex number $z$ with $|z|<1$, we have
 \begin{equation}\label{2.1}\sum_{k=1}^\infty kz^k=\f z{(1-z)^2}
 \ \ \t{and}\ \ \sum_{k=1}^\infty k^2z^k=\f{z(z+1)}{(1-z)^3}.
 \end{equation}
 \end{lemma}
 \Proof. This is easy. Recall the well-known identity
 $$\sum_{k=0}^\infty z^k=\f1{1-z}\quad (|z|<1)$$
 Taking derivatives of both sides, we get
 \begin{equation}\label{kz}\sum_{k=1}^\infty kz^{k-1}=\f1{(1-z)^2}
 \end{equation}
 and this implies the first identity in \eqref{2.1}.
 Taking derivatives of both sides of \eqref{kz}, we obtain
  \begin{equation}\label{k(k-1)z}\sum_{k=1}^\infty k(k-1)z^{k}=\f{2z^2}{(1-z)^3}.
 \end{equation}
 Adding this and the first identity in \eqref{2.1}, we immediately get the second identity
 in \eqref{2.1}. \qed
\medskip

 The beta function is connected with the Gamma function
 $$\Gamma(x):=\int_0^{+\infty}t^{x-1}e^{-t}dt\ \ (x>0)$$
 as first pointed out by Euler.

 \begin{lemma}[Euler]\label{Lem2.2}  For any $a>0$ and $b>0$, we have
\begin{equation} \label{B(a,b)} B(a,b)=\f{\Gamma(a)\Gamma(b)}{\Gamma(a+b)}.
\end{equation}
\end{lemma}

Now we present an auxiliary proposition.

\begin{proposition} Let $m> n>0$ be integers, and let $a,b,x$ be real numbers with $|x|<m^m/(n^n(m-n)^{m-n})$, and set
$$S_{m,n}(a,b,x):=\sum_{k=1}^\infty(ak+b)\f{x^k}{\bi{mk}{nk}}.$$
Then
$$S_{m,n}(a,b,x)=n\int_0^1 T_{m,n}(a,b,x;t)dt,$$
where
$$T_{m,n}(a,b,x;t):=t^{n-1}(1-t)^{m-n}x\f{(a-b)t^n(1-t)^{m-n}x+a+b}{(1-t^n(1-t)^{m-n}x)^3}.$$
\end{proposition}
\Proof. Clearly,
\begin{align*}S_{m,n}(a,b,x)&=\sum_{k=1}^\infty(ak+b)x^k\f{(nk)!((m-n)k)!}{(mk)!}
\\&=\sum_{k=1}^\infty(ak+b)x^k\f{nk\Gamma(nk)\Gamma((m-n)k+1)}{\Gamma(mk+1)}
\\&=n\sum_{k=1}^\infty(ak^2+bk)x^kB(nk,(m-n)k+1)
\\&=n\sum_{k=1}^\infty(ak^2+bk)x^k\int_0^1t^{nk-1}(1-t)^{(m-n)k}dt
\\&=n\int_0^1\f1t\sum_{k=1}^\infty(ak^2+bk)(t^n(1-t)^{(m-n)}x)^kdt.
\end{align*}
Note that  for $0\ls t\ls 1$ we have
$$\root m\of{\l(\f tn\r)^n\l(\f{1-t}{m-n}\r)^{m-n}}\ls \f {n\times\f tn+(m-n)\times\f {1-t}{m-n}}m=\f1m$$
and hence
$$|t^n(1-t)^{m-n}x|\ls\f{n^n(m-n)^{m-n}}{m^m}|x|<1.$$

Combining the above with Lemma \ref{Lem2.1}, we get
\begin{align*}&\f{S_{m,n}(a,b,x)}n
\\=\ &\int_0^1\l(a\f{t^{n-1}(1-t)^{m-n}x(t^n(1-t)^{m-n}x+1)}{(1-t^n(1-t)^{m-n}x)^3}
+b\f{t^{n-1}(1-t)^{m-n}x}{(1-t^n(1-t)^{m-n}x)^2}\r)dt
\\=\ &\int_0^1 T_{m,n}(a,b,x;t)dt.
\end{align*}
This concludes the proof. \qed

\section{Proof of Theorem \ref{C(3k,k)}}
 \setcounter{equation}{0}
 \setcounter{conjecture}{0}
 \setcounter{theorem}{0}
 \setcounter{proposition}{0}

 \begin{lemma} \label{Lem-arctan} For $-3<x<1$, we have
 \begin{equation}\label{q(x)}\arctan\f{x-1}{\sqrt{(1-x)(3+x)}}+\arctan\f{x+1}{\sqrt{(1-x)(3+x)}}
 =q(x),\end{equation}
 where $q(x)$ is as in Theorem \ref{C(3k,k)}.
 \end{lemma}
 \Proof. Let
 $$\al=\arctan\f{x-1}{\sqrt{(1-x)(3+x)}}\ \t{and}\ \beta=\arctan\f{x+1}{\sqrt{(1-x)(3+x)}}.$$
 If $x\not=-2$, then
 \begin{align*}\tan(\al+\beta)&=\f{\tan\al+\tan\beta}{1-(\tan\al)\tan\beta}
 \\&=\f{2x}{\sqrt{(1-x)(3+x)}}\l(1-\f{x^2-1}{(1-x)(3+x)}\r)^{-1}=\f x{x+2}\sqrt{\f{3+x}{1-x}}
 \end{align*}
 and hence $\al+\beta-\gamma\in\pi\Z,$
 where
 $$\gamma=\arctan\f x{x+2}\sqrt{\f{3+x}{1-x}}\in\l(-\f{\pi}2,\f{\pi}2\r).$$
 As $x<1$ we have $\al+\beta<\beta<\pi/2$ and hence $\al+\beta\ls\gamma$.
 If $-1\ls x<1$, then $\al+\beta\gs\al>-\pi/2$ and hence $\al+\beta=\gamma$.
 If $x\in(-2,-1)$, then $\gamma-\pi<-\pi<\al+\beta$ and hence $\al+\beta=\gamma$.

 Since
 $$\lim_{x\to-2\atop x>-2}\tan(\al+\beta)=\lim_{x\to-2\atop x>-2}\f x{x+2}\sqrt{\f{3+x}{1-x}}=-\infty,$$
 we have $\al+\beta=-\pi/2$ in the case $x=-2$.

 Now we consider the case $x\in(-3,-2)$. Note that $\al<0$ and $\beta<0$, but $\gamma>0$. So $\al+\beta=\gamma-\pi$.

 Combining the above, we obtain the desired identity \eqref{q(x)}. \qed

 \medskip
 \noindent{\it Proof of Theorem \ref{C(3k,k)}}. (i) Let
 \begin{align*}f(t)=&\ \f{(x+3)(2(t-1)(3t-2)x^3-9x+9)}{(1-x+t(1-t)^2x^3)^2}
 \\&\ +\f{2((t-1)(2t-1)x^4+(t-1)x^3-9x^2-9x+18)}{(x-1)(1-x+t(1-t)^2x^3)}
 \end{align*}
 and
 $$g(t)=\f{8x^2}{(x-1)\sqrt{(1-x)(3+x)}}\arctan\f{(2t-1)x-1}{\sqrt{(1-x)(3+x)}}.$$
 It is easy to verify that
 $$\f d{dt}\l(\f{(x-1)^2}{(x+3)^2}(f(t)+g(t))\r)=T_{3,1}\l((2x-3)^2,2x^2+2x-3,\f{x^3}{x-1},t\r).$$
 Thus, with the aid of Proposition 2.1, we get
 \begin{align*}&S_{3,1}\l((2x-3)^2,2x^2+2x-3,\f{x^3}{x-1}\r)
 \\=&\ \f{(x-1)^2}{(x+3)^2}(f(t)+g(t))\bigg|_{t=0}^1=\f{(x-1)^2}{(x+3)^2}(f(1)-f(0)+g(1)-g(0)).
 \end{align*}
 Note that
 \begin{align*}f(1)-f(0)&=\f{x+3}{(1-x)^2}(-2\times2x^3)+\f{2(-x^4+x^3)}{(x-1)(1-x)}
 \\&=-4x^3\f{x+3}{(1-x)^2}+\f{2x^3}{x-1}=-2x^3\f{x+7}{(1-x)^2}
 \end{align*}
 and
 \begin{align*}g(1)-g(0)=&\ \f{8x^2}{(x-1)\sqrt{(1-x)(3+x)}}
 \\&\ \times\l(\arctan\f{x-1}{\sqrt{(1-x)(3+x)}}+\arctan\f{x+1}{\sqrt{(1-x)(3+x)}}\r)
 \\=&\ \f{8x^2q(x)}{(x-1)\sqrt{(1-x)(3+x)}}
 \end{align*}
 with the help of Lemma \ref{Lem-arctan}. Therefore
 \begin{align*}&\sum_{k=1}^\infty\f{((2x-3)^2k+2x^2+2x-3)x^{3k}}{(x-1)^k\bi{3k}k}
 \\=\ &\f{(x-1)^2}{(x+3)^2}\l(-2x^3\f{x+7}{(x-1)^2}+\f{8x^2q(x)}{(x-1)\sqrt{(1-x)(3+x)}}\r)
 \\=\ &-2x^3\f{x+7}{(x+3)^2}+\f{8x^2(x-1)q(x)}{(x+3)^2\sqrt{(1-x)(3+x)}}.
 \end{align*}
 This proves Theorem 1.1(i).

 (ii) Set
 $$f_1(t)=\f{(x^3-13x^2+21x-9)(2x^3(3t^2-5t+2)-9x+9)}{(1-x+t(1-t)^2x^3)^2},$$
 \begin{align*}f_2(t)&=\f{2((2t^2-3t+1)x^5+6t(1-t)x^4+3(t-4)x^3+90x^2-135x+54)}{1-x+t(1-t)^2x^3}
 \end{align*}
 and
 $$f_3(t)=4x^2\l(2\log(1+(t-1)x)-\log(1+t^2x^2-tx(1+x))\r).$$
 It is easy to verify that
 $$\f d{dt}(x-1)(f_1(t)+f_2(t)+f_3(t))=T_{3,1}\l(s(x),t(x),\f{x^3}{x-1},t\r).$$
 Thus, by applying Proposition 2.1, we obtain
 \begin{align*}&S_{3,1}\l(s(x),t(x),\f{x^3}{x-1}\r)
 \\=\ &(x-1)(f_1(1)-f_1(0)+f_2(1)-f_2(0)+f_3(1)-f_3(0))
 \\=\ &(x-1)\l(\f{(x^3-13x^2+21x-9)(0-4x^3)}{(1-x)^2}+\f2{1-x}(-9x^3-(x^5-12x^3))\r)
 \\\ \ &+4x^2(x-1)(-\log(1-x)-2\log(1-x))
 \\=\ &-2x^3(x^2-24x+21)+12x^2(1-x)\log(1-x).
 \end{align*}
 Therefore
 \begin{align*}&\sum_{k=0}^\infty\f{(s(x)k+t(x))x^{3k}}{(x-1)^k\bi{3k}k}
 \\=&\ t(x)+S_{3,1}\l(s(x),t(x),\f{x^3}{x-1}\r)
 \\=&\ 2x^5-48x^4+69x^3-189x^2+243x-81
 \\\ &\ -2x^3(x^2-24x+21)+12x^2(1-x)\log(1-x).
 \\=&\ 27(x-1)(x^2-6x+3)+12x^2(1-x)\log(1-x).
 \end{align*}
 This proves the identity \eqref{log}.

 In view of the above, we have completed the proof of Theorem \ref{C(3k,k)}.
 \qed

 \section{Proof of Theorem \ref{Th1.1}}
 \setcounter{equation}{0}
 \setcounter{conjecture}{0}
 \setcounter{theorem}{0}
 \setcounter{proposition}{0}

 \begin{lemma} \label{Lem-x^(2k)} For $x>1/4$ we have
 \begin{equation}\label{arccoth}\sum_{k=0}^\infty\f1{x^{2k}\bi{4k}{2k}}
=\f{16x^2}{16x^2-1}+2x\l(\f{\mathrm{arccot}\sqrt{4x-1}}{(4x-1)\sqrt{4x-1}}
-\f{\mathrm{arccoth}\sqrt{4x+1}}{(4x+1)\sqrt{4x+1}}\r).
\end{equation}
 \end{lemma}
 \Proof. By Proposition 2.1,
 \begin{align*}\sum_{k=1}^\infty\f1{x^{2k}\bi{4k}{2k}}
 &=2\int_0^1T_{4,2}\l(0,1,\f1{x^2};t\r)dt
 \\&=\int_0^1\f{t(1-t)(1-2t)+t(1-t)}{x^2(1-t^2(1-t)^2/x^2)^2}dt
 \\&=\f1{2x^2(1-t^2(1-t)^2/x^2)}\bigg|_{t=0}^1+\f1{x^2}\int_0^1\f{t(1-t)}{(1-t^2(1-t)^2/x^2)^2}dt
 \\&=\f1{x^2}\int_0^{1/2}\f{t(1-t)}{(1-t^2(1-t)^2/x^2)^2}dt
 +\f1{x^2}\int_{1/2}^1\f{t(1-t)}{(1-t^2(1-t)^2/x^2)^2}dt
 \\&=\f{2}{x^2}\int_0^{1/2}\f{t(1-t)}{(1-t^2(1-t)^2/x^2)^2}dt.
 \end{align*}
 For $t\in[0,1/2]$, if we set $u=t(1-t)$ then
 $$t=\f{1-\sqrt{1-4u}}2\ \ \t{and}\ \ dt=\f{du}{\sqrt{1-4u}}.$$
 Thus
 \begin{equation}\label{1-4u}\sum_{k=1}^\infty\f1{x^{2k}\bi{4k}{2k}}
 =\f{2}{x^2}\int_0^{1/4}\f u{(1-u^2/x^2)^2\sqrt{1-4u}}du.
 \end{equation}

 Let $\psi(u)$ denote the expression
 $$\f{x(1+4u)\sqrt{1+4u}}{u^2-x^2}-\f{2(4x+1)}{\sqrt{4x-1}}\arctan\f{\sqrt{1-4u}}{\sqrt{4x-1}}
 +\f{2(4x-1)}{\sqrt{4x+1}}\mathrm{arctanh}\,\f{\sqrt{1-4u}}{\sqrt{4x+1}}.$$
 It is easy to verify that
 $$\f d{du}\l(\f x{16x^2-1}\psi(u)\r)=\f{2}{x^2}\cdot\f u{(1-u^2/x^2)^2\sqrt{1-4u}}.$$
 Combining this with \eqref{1-4u}, we get
 \begin{align*}\sum_{k=1}^\infty\f1{x^{2k}\bi{4k}{2k}}
 &=\f x{16x^2-1}\l(\psi\l(\f14\r)-\psi(0)\r)=-\f x{16x^2-1}\psi(0)
 \\&=-\f x{16x^2-1}\l(\f x{-x^2}-\f{2(4x+1)}{\sqrt{4x-1}}\arctan\f1{\sqrt{4x-1}}\r)
 \\&\ \ -\f x{16x^2-1}\cdot\f{2(4x-1)}{\sqrt{4x+1}}\mathrm{arctanh}\f{1}{\sqrt{4x+1}}
 \\&=\f1{16x^2-1}+\f{2x\ \mathrm{arccot}\sqrt{4x-1}}{(4x-1)\sqrt{4x-1}}
 -\f{2x\ \mathrm{arccoth}\sqrt{4x+1}}{(4x+1)\sqrt{4x+1}}
 \end{align*}
 and hence \eqref{arccoth} follows immediately. \qed

 \begin{remark} We can prove Lemma \ref{Lem-x^(2k)} in another way by noting that
 $$\sum_{k=0}^\infty\f1{x^{2k}\bi{4k}{2k}}=\f12\sum_{k=0}^\infty\f{1+(-1)^k}{x^k\bi{2k}k}
 \ \ \t{for}\ |x|>\f14,$$
 and using the identities \eqref{+arcsin} and \eqref{-arcsin} for $|x|\ls2$, which can be
 proved via Proposition 2.1.
 \end{remark}

 By Lemma \ref{Lem-x^(2k)},
 \begin{align*}\sum_{k=1}^\infty\f1{\bi{4k}{2k}}&=\f{45+25\pi\sqrt3-54\sqrt5\ \mathrm{arctanh}(1/\sqrt5)}{675},
\\\sum_{k=1}^\infty\f{4^k}{\bi{4k}{2k}}&=\f{12+9\pi-4\sqrt3\ \mathrm{arctanh}(1/\sqrt3)}{36},
\\\sum_{k=1}^\infty\f{9^k}{\bi{4k}{2k}}&=\f{189+98\pi\sqrt3-6\sqrt{21}\ \mathrm{arctanh}\sqrt{3/7}}{147},
\\
\sum_{k=1}^\infty\f{(9/4)^k}{\bi{4k}{2k}}&=\f 9{55}+\f{12}{55}\l(\f{11}{\sqrt{15}}\arctan\sqrt{\f 35}
-\f 5{\sqrt{33}}\ \mathrm{arctanh}\sqrt{\f3{11}}\r).
\end{align*}
We are also able to prove that
$$\sum_{k=1}^\infty\f1{k\bi{4k}{2k}}=\f{\sqrt3}9\pi-\f{2}5\sqrt5\log\f{1+\sqrt5}2$$
and
$$\sum_{k=1}^\infty\f1{k4^k\bi{4k}{2k}}=\f2{\sqrt7}\arctan\f1{\sqrt7}-\f{\log2}3$$
via the beta function.

 \medskip
 \noindent{\it Proof of Theorem \ref{Th1.1}}.
Taking derivatives of both sides of \eqref{arccoth}, we deduce that
\begin{equation}\label{k-arccoth}\begin{aligned}&\sum_{k=0}^\infty\f k{x^{2k}\bi{4k}{2k}}-\f{24x^2}{(16x^2-1)^2}
\\=\ &\f{x(1+2x)}{(4x-1)^2\sqrt{4x-1}}\mathrm{arccot}\sqrt{4x-1}
+\f{x(1-2x)}{(4x+1)^2\sqrt{4x+1}}\mathrm{arccoth}\sqrt{4x+1}.
\end{aligned}\end{equation}
Via $2(4x+1)\times$\eqref{k-arccoth}$+(1-2x)\times$\eqref{arccoth} we see that \eqref{comput-arccot} holds.
Similarly, via $2(4x-1)\times$\eqref{k-arccoth}$-(2x+1)\times$\eqref{arccoth} we obtain
$$\sum_{k=0}^\infty\f{2(4x-1)k-2x-1}{x^{2k}\bi{4k}{2k}}
=\f{8x^2}{(4x+1)^2}\l(\f3{\sqrt{4x+1}}\mathrm{arccoth}\sqrt{4x+1}-4x-4\r),$$
which is equivalent to \eqref{comput-arccoth} since
\begin{align*}\f{\mathrm{arccoth}\sqrt{4x+1}}{\sqrt{4x+1}}
&=\sum_{k=0}^\infty\f{1}{(2k+1)(\sqrt{4x+1})^{2k+2}}
\\&=\f1{4x+1}\sum_{k=0}^\infty\f1{(2k+1)(4x+1)^k}=\f{R(4x+1)}{4x+1}.
\end{align*}

Putting $x=1,1/2,1/3,2/3$ in \eqref{comput-arccot} we immediately get \eqref{new-1}-\eqref{new-9/4}.
In light of \eqref{R(x)}, the identity \eqref{New3} follows from \eqref{comput-arccoth} with $x=2$.

In view of the above, we have completed the proof of Theorem \ref{Th1.1}. \qed

\section{Proof of Theorem \ref{C(4k,2k)}}
 \setcounter{equation}{0}
 \setcounter{conjecture}{0}
 \setcounter{theorem}{0}
 \setcounter{proposition}{0}

 \begin{lemma} \label{Lem-t} For any $u<1$ with $u\not=0$, we have
 \begin{equation}\label{u}
\sum_{k=0}^\infty \f{u^k}{2k+1}\bg(\bg(1-i\sqrt{1-u}\bg)^{-2k-1}+\bg(1+i\sqrt{1-u}\bg)^{-2k-1}\bg)
=\f{\mathrm{arctanh} \sqrt{u}}{\sqrt u}.
\end{equation}
 \end{lemma}
 \Proof. It suffices to prove that
 \begin{equation}\label{t-series}
\sum_{k=0}^\infty \f{t^{2k+1}}{2k+1}\bg(\bg(1-i\sqrt{1-t^2}\bg)^{-2k-1}+\bg(1+i\sqrt{1-t^2}\bg)^{-2k-1}\bg)
=\mathrm{arctanh}\,t
\end{equation}
 for each $t\in\C$ with $t^2<1$.
 Note that
$$\bg|\f t{1\pm i\sqrt{1-t^2}}\bg|^2=\f{|(1-t^2)-1|}{1+(1-t^2)}<1.$$

Let $f(t)$ and $g(t)$ denote the left-hand side and the right-hand side of \eqref{t-series} respectively.
Then
\begin{align*}f'(t)&=\sum_{k=0}^\infty t^{2k}\bg(\bg(1-i\sqrt{1-t^2}\bg)^{-2k-1}+\bg(1+i\sqrt{1-t^2}\bg)^{-2k-1}\bg)
\\&\ \ -\sum_{k=0}^\infty t^{2k+1}\bg(\bg(1-i\sqrt{1-t^2}\bg)^{-2k-2}-\bg(1+i\sqrt{1-t^2}\bg)^{-2k-2}
\bg)\f{it}{\sqrt{1-t^2}}
\\&=\f1{1-i\sqrt{1-t^2}}\sum_{k=0}^\infty\l(\f{t^2}{(1-i\sqrt{1-t^2})^2}\r)^k
\\&\ \ +\f1{1+i\sqrt{1-t^2}}\sum_{k=0}^\infty\l(\f{t^2}{(1+i\sqrt{1-t^2})^2}\r)^k
\\&\ \ -\f i{\sqrt{1-t^2}}\sum_{k=0}^\infty\l(\l(\f{t^2}{(1-i\sqrt{1-t^2})^2}\r)^{k+1}
-\l(\f{t^2}{(1+i\sqrt{1-t^2})^2}\r)^{k+1}\r).
\end{align*}
Note that $(1\pm i\sqrt{1-t^2})^2=t^2\pm 2i\sqrt{1-t^2}$. Therefore
\begin{align*}f'(t)&=\f1{1-i\sqrt{1-t^2}}\cdot\f1{1-t^2/(t^2-2i\sqrt{1-t^2})}
\\&\ \ +\f1{1+i\sqrt{1-t^2}}\cdot\f1{1-t^2/(t^2+2i\sqrt{1-t^2})}
\\&\ \ -\f i{\sqrt{1-t^2}}\cdot\f{t^2}{t^2-2i\sqrt{1-t^2}}\cdot\f1{1-t^2/(t^2-2i\sqrt{1-t^2})}
\\&\ \ +\f i{\sqrt{1-t^2}}\cdot\f{t^2}{t^2+2i\sqrt{1-t^2}}\cdot\f1{1-t^2/(t^2+2i\sqrt{1-t^2})}.
\end{align*}
Hence
\begin{align*}f'(t)&=\ \f{1-i\sqrt{1-t^2}}{-2i\sqrt{1-t^2}}+\f{1+i\sqrt{1-t^2}}{2i\sqrt{1-t^2}}
\\&\ \ -\f{it^2}{\sqrt{1-t^2}}\cdot\f1{-2i\sqrt{1-t^2}}+\f{it^2}{\sqrt{1-t^2}}\cdot\f1{2i\sqrt{1-t^2}}
\\&=1+\f{t^2}{1-t^2}=\f1{1-t^2}=g'(t).
\end{align*}
Thus $f(t)-g(t)$ is a constant. Since $f(0)=0=g(0)$, we have $f(t)=g(t)$.
This concludes our proof of Lemma \ref{Lem-t}. \qed

 \begin{lemma}\label{Lem-FG} Let $x>1$ or $x<0$, and let
 \begin{equation}
 \label{F(s)} F(s)=\f{20x^2-56x+35}{(4s^2+x(x-1))^2}-\f{4(x-1)(2x-3)(2x-1)^2}{(4s^2+x(x-1))^3}.
 \end{equation}
 Then
 \begin{equation}\label{sF(s)}
 8x^2(1-x)\int_{0}^{1/4}\f{sF(s)}{\sqrt{1-4s}}ds=3(x-1)R(x)+5x-6.
 \end{equation}
 \end{lemma}
 \Proof. Note that $x(x-1)>0$.
 Let $G(s)$ denote the expression
 \begin{align*}\ &(24s^3+4s^2x+2sx(x-1)(8x-9)+x(x-1)(5x-6))\f{\sqrt{1-4s}\sqrt x}{(4s^2+x(1-x))^2}
 \\\ &+\f{3((1-x)i+\sqrt{x}\sqrt{x-1})}{\sqrt{x-1}\sqrt{1-2i\sqrt{x}\sqrt{x-1}}}\ \mathrm{arctanh}\l( \f{\sqrt{1-4s}}{\sqrt{1-2i\sqrt{x}\sqrt{x-1}}}\r)
 \\\ &+\f{3((x-1)i+\sqrt{x}\sqrt{x-1})}{\sqrt{x-1}\sqrt{1+2i\sqrt{x}\sqrt{x-1}}}\ \mathrm{arctanh} \l(\f{\sqrt{1-4s}}{\sqrt{1+2i\sqrt{x}\sqrt{x-1}}}\r).
 \end{align*}
As $$\f d{dz}(\mathrm{arctanh}\, z)=\f1{1-z^2},$$
it is routine to verify that
$$\f d{ds}\l(\f{G(s)}{8x^{3/2}}\r)=\f{sF(s)}{\sqrt{1-4s}}.$$
Actually we find the expression of $G(s)$ by {\tt Mathematica}.
Since $x>1$ or $x<0$, we have
$$\l|\sqrt{1\pm2i\sqrt x\sqrt{x-1}}\r|=\l|\sqrt x\pm i\sqrt{x-1}\r|=\sqrt{|2x-1|}>1.$$
Thus
\begin{align*}&8x^{3/2}\int_0^{1/4}\f{sF(s)}{\sqrt{1-4s}}ds
\\=&\ G\l(\f14\r)-G(0)=-G(0)
\\=&\ -\f{x(x-1)(5x-6)\sqrt x}{x^2(1-x)^2}
\\&\ -\f{3((1-x)i+\sqrt{x}\sqrt{x-1})}{\sqrt{x-1}}
\sum_{k=0}^\infty\f1{(2k+1)\l(\sqrt{1-2i\sqrt{x}\sqrt{x-1}}\r)^{2k+2}}
\\&\ -\f{3((x-1)i+\sqrt{x}\sqrt{x-1})}{\sqrt{x-1}}
\sum_{k=0}^\infty\f1{(2k+1)\l(\sqrt{1+2i\sqrt{x}\sqrt{x-1}}\r)^{2k+2}}
\\=&\ \f{5x-6}{x(1-x)}\sqrt x
\\\ &-3\sum_{k=0}^\infty\l(\f{\sqrt{x}-i\sqrt{x-1}}
{(2k+1)(1-2i\sqrt{x}\sqrt{x-1})^{k+1}}+\f{\sqrt{x}+i\sqrt{x-1}}
{(2k+1)(1+2i\sqrt{x}\sqrt{x-1})^{k+1}}\r)
\end{align*}
which coincides with
\begin{align*}
&\ \f{5x-6}{x(1-x)}\sqrt x-3\sum_{k=0}^\infty\f1{2k+1}\l(\f{\sqrt{x}-i\sqrt{x-1}}
{(\sqrt x-i\sqrt{x-1})^{2k+2}}+\f{\sqrt{x}-i\sqrt{x-1}}
{(\sqrt x-i\sqrt{x-1})^{2k+2}}\r)
\\&\ \ = \f{5x-6}{x(1-x)}\sqrt x-6\sum_{k=0}^\infty\f{\sqrt x}{(2k+1)x^{k+1}}\Re\bg(\bg(1-i\sqrt{\f{x-1}x}\bg)^{-2k-1}\bg).
\end{align*}
This reduces the desired identity \eqref{sF(s)} to
\begin{equation*}
\sum_{k=0}^\infty \f{x^{-k}}{2k+1}\bg(\bg(1-i\sqrt{\f{x-1}x}\bg)^{-2k-1}+\bg(1+i\sqrt{\f{x-1}x}\bg)^{-2k-1}\bg)
=R(x)
\end{equation*}
which follows from the identity \eqref{u} with $u=1/x$. This ends our proof. \qed

\begin{lemma}\label{int-F(s)} Let $x>1$ or $x<0$, and let $F(s)$ be defined by \eqref{F(s)}. Then
$$\int_0^1t(1-t)^2F(t(1-t))dt=\int_0^{1/4}\f{sF(s)}{\sqrt{1-4s}}ds.$$
\end{lemma}
\Proof. Note that
$$\int_0^1\f{8t(1-t)(1-2t)}{(4t^2(1-t)^2+x(x-1))^2}dt
=\int_0^1 \l(\f{-1}{4t^2(1-t)^2+x(x-1)}\r)'dt=0$$
and
$$\int_0^1\f{8t(1-t)(1-2t)}{(4t^2(1-t)^2+x(x-1))^3}dt
=\int_0^1 \l(\f{-1/2}{(4t^2(1-t)^2+x(x-1))^2}\r)'dt=0.$$
Thus
\begin{align*}&2\int_0^1t(1-t)^2F(t(1-t))dt
\\=\ &\int_0^1t(1-t)((1-2t)+1)F(t(1-t))dt
\\=\ &\int_0^{1/2}t(1-t)F(t(1-t))dt+\int_{1/2}^1 u(1-u)F(u(1-u))du
\\=\ &\int_0^{1/2}t(1-t)F(t(1-t))dt+\int_{1/2}^0 (1-t)tF((1-t)t)d(1-t)
\\=\ &2\int_0^{1/2}t(1-t)F(t(1-t))dt.
\end{align*}
For $t\in[0,1/2]$, if we set $s=t(1-t)$, then $t=(1-\sqrt{1-4s})/2$ and hence
$$dt=-\f14\cdot\f{-4}{\sqrt{1-4s}}ds=\f{ds}{\sqrt{1-4s}}.$$
Therefore
$$\int_0^1t(1-t)^2F(t(1-t))dt=\int_0^{1/2}t(1-t)F(t(1-t))dt=\int_0^{1/4}\f{sF(s)}{\sqrt{1-4s}}ds$$
as desired. \qed

 \medskip
 \noindent{\it Proof of Theorem \ref{C(4k,2k)}}.
 Note that
 $$\bi{4k}{2k}\sim \f{16^k}{\sqrt{2k\pi}}\ \ \t{and}\ \  0< \f4{x(x-1)}<16.$$
 So the series in \eqref{gen} converges absolutely.
 By Propsotion 2.1, we have
 \begin{align*}&\sum_{k=1}^\infty \f{(2(2x-1)^2(2x-3)k-(4x^3-16x^2+7x+6))4^k}{(x(1-x))^k\bi{4k}{2k}}
 \\=\ &2\int_0^1T_{4,2}\l(2(2x-1)^2(2x-3),-(4x^3-16x^2+7x+6),\f 4{x(1-x)},t\r)dt.
 \end{align*}
 It is easy to verify that
 \begin{align*}&T_{4,2}\l(2(2x-1)^2(2x-3),-(4x^3-16x^2+7x+6),\f 4{x(1-x)},t\r)
\\=\ &\f{4x^2(x-1)(2-x^2-56x+35)t(1-t)^2}{(4t^2(1-t)^2+x(x-1))^2}
\\ \ \ &-\f{16x^2(x-1)^2(2x-1)^2(2x-3)t(1-t)^2}{(4t^2(1-t)^2+x(x-1))^3}
\\=\ &4x^2(x-1)t(1-t)^2F(t(1-t)),
\end{align*}
where the function $F$ is given by \eqref{F(s)}.
Combining this with Lemmas \ref{Lem-FG} and \ref{int-F(s)}, we get
\begin{align*}&\sum_{k=1}^\infty \f{(2(2x-1)^2(2x-3)k-(4x^3-16x^2+7x+6))4^k}{(x(1-x))^k\bi{4k}{2k}}
\\=\ &8x^2(x-1)\int_0^1 t(1-t)^2F(t(1-t))dt=8x^2(x-1)\int_0^{1/4}\f{sF(s)}{\sqrt{1-4s}}ds
\\=\ &3(1-x)R(x)-5x+6
\end{align*}
and hence
\begin{align*}&\sum_{k=0}^\infty \f{(2(2x-1)^2(2x-3)k-(4x^3-16x^2+7x+6))4^k}{(x(1-x))^k\bi{4k}{2k}}
\\=\ &-(4x^3-16x^2+7x+6)+3(1-x)R(x)-5x+6
\\=\ &(1-x)(3R(x)+4x(x-3)).
\end{align*}
This proves \eqref{gen}.

As $x>(1+\sqrt2)/2$ or $x<(1-\sqrt2)/2$, we see that $1-x<(1-\sqrt2)/2$ or $1-x>(1+\sqrt2)/2$.
Note that \eqref{gen} with $x$ replaced by $1-x$ yields \eqref{gen-dual}.
This concludes our proof of Theorem \ref{C(4k,2k)}. \qed

\section{Conjectural series with summands containing $\bi{2k}k/(\bi{3k}k\bi{6k}{3k})$}
 \setcounter{equation}{0}
 \setcounter{conjecture}{0}
 \setcounter{theorem}{0}
 \setcounter{proposition}{0}

 In 2013, the author \cite{S13} proved that
 $$2(2k+1)\bi{2k}k\ \bigg|\ \bi{3k}k\bi{6k}{3k}\quad \ \t{for all}\ k\in\N.$$
 In 2014 W. Chu and W. Zhang \cite[Example 27]{CZ} obtained an identity which has the following equivalent form:
 $$\sum_{k=1}^\infty\f{(7k-1)(-4)^k\bi{2k}k}{(2k-1)k\bi{3k}k\bi{6k}{3k}}=-\f{\pi}4;$$
 its corresponding $p$-adic congruences were conjectured by Sun \cite[Conjecture 4.7]{BNK}.
 Motivated by this and the author's recent work \cite{harmonic},
 in this section we pose some conjectures on series whose summands involve $\bi{2k}k/(\bi{3k}k\bi{6k}{3k})\ (k\in\N)$.

 As usual, those rational numbers
 $$H_n:=\sum_{0<k\ls n}\f1k\ \ (n=0,1,2,\ldots)$$
 are called {\it harmonic numbers}. For $m=2,3,\ldots$, those rational numbers
 $$H_n^{(m)}:=\sum_{0<k\ls n}\f1{k^m}\ \ (n=0,1,2,\ldots)$$
 are called {\it harmonic numbers of order $m$}.

 {\it Dirichlet's beta function} is defined by
 $$\beta(s)=\sum_{k=0}^\infty\f{(-1)^k}{(2k+1)^s}\ \ (s=1,2,3,\ldots).$$
 It is well known that $\beta(1)=\pi/4$. Note that $G=\beta(2)$ is {\it Catalan's constant}.
 For a series $\sum_{k=0}^\infty a_k$ with $a_0,a_1,\ldots$ real numbers, if $\lim_{k\to+\infty}a_{k+1}/a_k=r\in(-1,1)$, then we say that its {\it converging rate} is $r$.

\begin{conjecture} [2023-08-23] \label{7k-1Conj} {\rm (i)} We have
\begin{align}\sum_{k=1}^\infty\f{(-4)^k\bi{2k}k((7k-1)H_{k-1}-(6k-1)/(4k-2))}{k(2k-1)\bi{3k}k\bi{6k}{3k}}
&=\pi\log2-2G,
\\\sum_{k=1}^\infty\f{(-4)^k\bi{2k}k((7k-1)H_{2k-1}-9(6k-1)/(8k-4))}{k(2k-1)\bi{3k}k\bi{6k}{3k}}
&=\f34\pi\log2-G,
\end{align}
and
\begin{equation}\begin{aligned}
&\sum_{k=1}^\infty\f{(-4)^k\bi{2k}k}{k(2k-1)\bi{3k}k\bi{6k}{3k}}\l((7k-1)(2H_{6k-1}-H_{3k-1})-\f{34k-9}{4k-2}\r)
\\&\qquad=\f{\pi}2\log2-2G.
\end{aligned}\end{equation}

{\rm (ii)} Let $p$ be an odd prime. Then
\begin{equation*}\begin{aligned}\sum_{k=0}^{(p-3)/2}\f{\bi{3k}k\bi{6k}{3k}((7k+1)H_k+(6k+1)/(4k+2))}
{(2k+1)(-4)^k\bi{2k}k}
\eq\l(\f{-1}p\r)q_p(2)\pmod{p}
\end{aligned}\end{equation*}
and
\begin{equation*}\begin{aligned} &\sum_{k=0}^{(p-3)/2}\f{\bi{3k}k\bi{6k}{3k}((7k+1)H_{2k}+9(6k+1)/(8k+4))}
{(2k+1)(-4)^k\bi{2k}k}
\\&\quad\eq\l(\f{-1}p\r)\f34\l(q_p(2)-p\,q_p(2)^2\r)\pmod{p^2},
\end{aligned}\end{equation*}
where $q_p(2)$ denotes the Fermat quotient $(2^{p-1}-1)/p$.
\end{conjecture}
\begin{remark} \label{Rem(-4)^k}
Suitable linear combinations of the three identities in Conjecture \ref{7k-1Conj}(i) yield the identities
\begin{equation}\begin{aligned}&\ \sum_{k=1}^\infty\f{(-4)^k\bi{2k}k((7k-1)(30H_{6k-1}-15H_{3k-1}-14H_{2k-1}+3H_{k-1})-75/2)}{(2k-1)k\bi{3k}k\bi{6k}{3k}}
\\&\qquad\qquad=-22G
\end{aligned}
\end{equation}
and
\begin{equation}\begin{aligned}&\ \sum_{k=1}^\infty\f{(-4)^k\bi{2k}k((7k-1)(16H_{6k-1}-8H_{3k-1}-6H_{2k-1}-5H_{k-1})-20)}{(2k-1)k\bi{3k}k\bi{6k}{3k}}
\\&\qquad\quad=-\f{11}2\pi\log2.
\end{aligned}
\end{equation}
\end{remark}

\begin{conjecture} [2023-10-28]\label{7k-1High} {\rm (i)} We have
\begin{equation}\sum_{k=1}^\infty\f{(-4)^k\bi{2k}k((7k-1)(6H_{2k-1}^{(2)}-H_{k-1}^{(2)})-9(6k-1)/(2k-1)^2)}
{(2k-1)k\bi{3k}k\bi{6k}{3k}}=\f{\pi^3}{24}
\end{equation}
and
\begin{equation}\begin{aligned}&\sum_{k=1}^\infty\f{(-4)^k\bi{2k}k((7k-1)(8H_{2k-1}^{(3)}-H_{k-1}^{(3)})
-12(6k-1)/(2k-1)^3)}
{(2k-1)k\bi{3k}k\bi{6k}{3k}}
\\&\qquad=\f 94\pi\zeta(3)-\f{48}7\beta(4).
\end{aligned}
\end{equation}

{\rm (ii)} For any prime $p>3$, we have
\begin{align*}&\sum_{k=0}^{(p-3)/2}\f{\bi{3k}k\bi{6k}{3k}}
{(2k+1)(-4)^k\bi{2k}k}\l((7k+1)(6H_{2k}^{(2)}-H_k^{(2)})+\f{9(6k+1)}{(2k+1)^2}\r)
\\&\qquad\eq E_{p-3}\pmod p,
\end{align*}
where $E_0,E_1,E_2,\ldots$ are the Euler numbers.
For each odd prime $p\not=7$, we have
\begin{align*}&\sum_{k=0}^{(p-3)/2}\f{\bi{3k}k\bi{6k}{3k}}
{(2k+1)(-4)^k\bi{2k}k}\l((7k+1)(8H_{2k}^{(3)}-H_k^{(3)})
+\f{12(6k+1)}{(2k+1)^3}\r)
\\&\qquad\eq \l(\f{-1}p\r)\f34B_{p-3}\pmod p,
\end{align*}
where $B_0,B_1,B_2,\ldots$ are the Bernoulli numbers.
\end{conjecture}
\begin{remark} The two series in Conjecture \ref{7k-1High}(i) have converging rate $-1/27$.
\end{remark}

\begin{conjecture} [2023-09-09] We have
\begin{equation}\sum_{k=1}^\infty\f{(-4)^k\bi{2k}k
((280k-51)H(k)-1352)}{k\bi{3k}k\bi{6k}{3k}}=18\pi-624G
\end{equation}
and
\begin{equation}\sum_{k=1}^\infty\f{(-4)^k\bi{2k}k
(17(952k-201)H(k)-50924)}{\bi{3k}k\bi{6k}{3k}}=-452-1587\pi-26520G.
\end{equation}
where
$$H(k):=30H_{6k-1}-15H_{3k-1}-22H_{2k-1}+9H_{k-1}.$$
\end{conjecture}
\begin{remark}\label{Rem} In view of the first identity in this section, we have
$$\sum_{k=1}^\infty\f{(280k-51)(-4)^k\bi{2k}k}{k\bi{3k}k\bi{6k}{3k}}=-6\pi-10,$$
because $(2k-1)(280k-51)-24(7k-1)=5(112k^2-110k+15)$, and
$$\sum_{k=1}^n\f{(112k^2-110k+15)(-4)^k\bi{2k}k}{(2k-1)k\bi{3k}k\bi{6k}{3k}}=-2+\f{(-1)^n2^{2n+1}\bi{2n}n}{\bi{3n}n\bi{6n}{3n}}$$
tends to $-2$ as $n\to+\infty$. Similarly, we have
$$\sum_{k=1}^\infty\f{(952k-201)(-4)^k\bi{2k}k}{\bi{3k}k\bi{6k}{3k}}=-15\pi-42,$$
because $2k(952k-201)-5(280k-51)=17(112k^2-106k+15)$, and
$$\sum_{k=1}^n\f{(112k^2-106k+15)(-4)^k\bi{2k}k}{k\bi{3k}k\bi{6k}{3k}}=-2+\f{(-1)^n(2n+1)2^{2n+1}\bi{2n}n}{\bi{3n}n\bi{6n}{3n}}$$
tends to $-2$ as $n\to+\infty$.
\end{remark}

\begin{conjecture} [2023-09-09] {\rm (i)} We have
\begin{equation}\label{8^k/}\sum_{k=0}^\infty\f{(350k-17)8^k\bi{2k}k}
{\bi{3k}k\bi{6k}{3k}}=15\sqrt2\pi+27.
\end{equation}
Also,
\begin{equation}\begin{aligned}&\sum_{k=1}^\infty\f{8^k\bi{2k}k}{\bi{3k}k\bi{6k}{3k}}
\l(21(350k-17)(2H_{6k-1}-H_{3k-1}-H_{k-1})+4850\r)
\\&\qquad=976+1020\sqrt2\pi+945\pi\sqrt2\log2
\end{aligned}
\end{equation}
and
\begin{equation}\begin{aligned}&\sum_{k=1}^\infty\f{8^k\bi{2k}k
}{\bi{3k}k\bi{6k}{3k}}\l(7(350k-17)(H_{2k-1}-H_{k-1})+2225\r)
\\&\ = 276+\f{493}{\sqrt2}\pi+\f{315}{\sqrt2}\pi\log2-420L,
\end{aligned}\end{equation}
where
$$L:=L\l(2,\l(\f{-8}{\cdot}\r)\r)=\sum_{n=1}^\infty\f{(\f{-8}n)}{n^2}
=\sum_{k=0}^\infty\f{(-1)^{k(k-1)/2}}{(2k+1)^2}.$$

{\rm (ii)} Let $p$ be any odd prime. Then
\begin{align*}\sum_{k=1}^{(p-1)/2}\f{k(350k+17)\bi{3k}k\bi{6k}{3k}}{(2k+1)8^k\bi{2k}k}&\eq15\l(\f{-2}p\r)
-\l(\f 2p\r)\f{93}2p\pmod{p^2}.
\end{align*}
\end{conjecture}
\begin{remark} The series in \eqref{8^k/} has converging rate $2/27$.
\end{remark}

\begin{conjecture} [2023-09-23] {\rm (i)} We have
\begin{align}&\sum_{k=1}^\infty\f{\bi{2k}k8^k((50k-7)(H_{2k-1}-H_{k-1})+5)}{k\bi{3k}k\bi{6k}{3k}}
 = 3\sqrt2\pi(1+\log2)-8L.
\end{align}
and
\begin{equation}\begin{aligned}&\sum_{k=1}^\infty\f{\bi{2k}k8^k((50k-7)(2H_{6k-1}-H_{3k-1}-H_{k-1})-10)}{k\bi{3k}k\bi{6k}{3k}}
=\sqrt2\pi(4+6\log2).
\end{aligned}
\end{equation}

{\rm (ii)} We have
\begin{equation}\begin{aligned}&\sum_{k=1}^\infty\f{\bi{2k}k8^k((5k-1)(H_{2k-1}-H_{k-1})-3(6k-1)/(8k-4))}{(2k-1)k\bi{3k}k\bi{6k}{3k}}
\\&\qquad=\f38\sqrt2\,\pi\log2-L,
\end{aligned}
\end{equation}
\begin{equation}\begin{aligned}&\sum_{k=1}^\infty\f{\bi{2k}k8^k((5k-1)(2H_{6k-1}-H_{3k-1}-H_{k-1})-(6k-2)/(2k-1))}{(2k-1)k\bi{3k}k\bi{6k}{3k}}
\\&\qquad=\f34\sqrt2\,\pi\log2,
\end{aligned}
\end{equation}
and
\begin{equation}\begin{aligned}&\sum_{k=1}^\infty\f{\bi{2k}k8^k((5k-1)(12H_{6k-1}-6H_{3k-1}-4H_{2k-1}-2H_{k-1})-9)}{(2k-1)k\bi{3k}k\bi{6k}{3k}}
\\&\qquad=3\sqrt2\,\pi\log2+4L.
\end{aligned}
\end{equation}

{\rm (iii)} We have the identity
\begin{equation} \sum_{k=1}^\infty\f{\bi{2k}k8^k((5k-1)(16H_{2k-1}^{(2)}-3H_{k-1}^{(2)})-12(6k-1)/(2k-1)^2)}
{(2k-1)k\bi{3k}k\bi{6k}{3k}}=\f{\pi^3}{12\sqrt2}.
\end{equation}
Also,
\begin{align*}&\sum_{k=0}^{(p-3)/2}\f{\bi{3k}k\bi{6k}{3k}((5k+1)(16H_{2k}^{(2)}-3H_k^{(2)})+12(6k+1)/(2k+1)^2)}
{(2k+1)\bi{2k}k8^k}
\\&\qquad\eq-\f14E_{p-3}\l(\f14\r)\pmod p
\end{align*}
for any prime $p>3$,
where $E_{p-3}(x)$ denotes the Euler polynomial of degree $p-3$.
\end{conjecture}
\begin{remark} In the spirit of the arguments in Remark \ref{Rem}, \eqref{8^k/} implies that
\begin{equation}\sum_{k=1}^\infty\f{(50k-7)8^k\bi{2k}k}{k\bi{3k}k\bi{6k}{3k}}=4+2\sqrt2\pi
\ \t{and}\ \sum_{k=1}^\infty\f{(5k-1)8^k\bi{2k}k}{k(2k-1)\bi{3k}k\bi{6k}{3k}}=\f{\pi}{2\sqrt2}.
\end{equation}
\end{remark}

\begin{conjecture}[2023-10-18]  We have
\begin{equation}\label{-3^k}\sum_{k=1}^\infty\f{(130k-21)\bi{2k}k}{k(2k-1)(-3)^k\bi{3k}k\bi{6k}{3k}}
=-\f{2\pi}{3\sqrt3}.
\end{equation} Also,
\begin{equation}\begin{aligned}&\sum_{k=1}^\infty\f{\bi{2k}k\l((130k-21)(H_{2k-1}-H_{k-1})-26(6k-1)/(2k-1)\r)}{k(2k-1)(-3)^k\bi{3k}k\bi{6k}{3k}}
\\&\qquad=2K-\f{2\pi}{3\sqrt3}\log3\end{aligned}
\end{equation}
and
\begin{equation}\begin{aligned}&\sum_{k=1}^\infty\f{\bi{2k}k\l((130k-21)(2H_{6k-1}-H_{3k-1}-H_{k-1})-16(13k-4)/(2k-1)\r)}{k(2k-1)(-3)^k\bi{3k}k\bi{6k}{3k}}
\\&\qquad=K+\f{2\pi}{3\sqrt3}\log3,\end{aligned}
\end{equation}
where
$$K:=L\l(2,\l(\f{-3}{\cdot}\r)\r)=\sum_{n=1}^\infty\f{(\f n3)}{n^2}
=\sum_{k=0}^\infty\l(\f1{(3k+1)^2}-\f1{(3k+2)^2}\r).$$
\end{conjecture}
\begin{remark} The series in \eqref{-3^k} has converging rate $-1/324$.
A linear combination of the last two formulae yields the identity
\begin{align*}&\sum_{k=1}^\infty\f{\bi{2k}k((130k-21)(26H_{6k-1}-13H_{3k-1}-10H_{2k-1}-3H_{k-1})-572)}
{k(2k-1)(-3)^k\bi{3k}k\bi{6k}{3k}}
\\&\qquad=-33K-\f{2\pi}{\sqrt3}\log3.
\end{align*}
\end{remark}

\begin{conjecture}[2023-10-18]  We have
\begin{equation}\label{-27^k}\sum_{k=1}^\infty\f{(10k-1)(-27)^k\bi{2k}k}{k(2k-1)\bi{3k}k\bi{6k}{3k}}
=-\f{4\pi}{\sqrt3}.
\end{equation} Also,
\begin{equation}\sum_{k=1}^\infty\f{\bi{2k}k(-27)^k\l((10k-1)H_{k-1}+2(6k-1)/(6k-3)\r)}{k(2k-1)\bi{3k}k\bi{6k}{3k}}
=2\pi\sqrt3\log3-18K,
\end{equation}
\begin{equation}\sum_{k=1}^\infty\f{\bi{2k}k(-27)^k\l((10k-1)H_{2k-1}-8(6k-1)/(6k-3)\r)}{k(2k-1)\bi{3k}k\bi{6k}{3k}}
=2\pi\sqrt3\log3-9K,
\end{equation}
and
\begin{equation}\begin{aligned}&\sum_{k=1}^\infty\f{\bi{2k}k(-27)^k\l((10k-1)(2H_{6k-1}-H_{3k-1})
-2(102k-29)/(18k-9)\r)}{k(2k-1)\bi{3k}k\bi{6k}{3k}}
\\&\qquad=15K-\f2{\sqrt3}\pi\log3.\end{aligned}
\end{equation}
\end{conjecture}
\begin{remark} The series in \eqref{-27^k} has converging rate $-1/4$.
Suitable combinations of the last three formulae yield that
$$\sum_{k=1}^\infty\f{\bi{2k}k(-27)^k(10k-1)(H_{2k-1}+4H_{k-1})}{k(2k-1)\bi{3k}k\bi{6k}{3k}}
=10\pi\sqrt3\log3-81K$$
and
\begin{align*}&\sum_{k=1}^\infty\f{\bi{2k}k(-27)^k((10k-1)(6H_{6k-1}-3H_{3k-1}+11H_{k-1})-12)}
{k(2k-1)\bi{3k}k\bi{6k}{3k}}
\\&\qquad=24\pi\sqrt3\log3-243K.
\end{align*}
\end{remark}

\begin{conjecture}[2023-09-28] We have
\begin{equation}\sum_{k=1}^\infty\f{\bi{2k}k16^k(46k^2-11k+1)}{k^2(2k-1)^2\bi{3k}k\bi{6k}{3k}}=2\pi^2,
\end{equation}
\begin{equation}\label{H2k-1}\sum_{k=1}^\infty\f{\bi{2k}k16^k((46k^2-11k+1)H_{2k-1}+5k(6k-1)/(2k-1))}{k^2(2k-1)^2\bi{3k}k\bi{6k}{3k}}
=28\zeta(3),
\end{equation}
and
\begin{equation}\label{H631}\begin{aligned}&
\sum_{k=1}^\infty\f{\bi{2k}k16^k\l((46k^2-11k+1)(2H_{6k-1}-H_{3k-1}+H_{k-1})+\f{8k(19k-2)}{2k-1}\r)}{k^2(2k-1)^2\bi{3k}k\bi{6k}{3k}}
\\&\qquad=112\zeta(3).\end{aligned}\end{equation}
Also,
\begin{equation}\begin{aligned}&\sum_{k=1}^\infty\f{\bi{2k}k16^k((46k^2-11k+1)(21H_{2k-1}-5H_{k-1})+5(2k-1)^2/(2k))}{k^2(2k-1)^2\bi{3k}k\bi{6k}{3k}}
\\&\qquad =238\zeta(3)+20\pi^2\log2
\end{aligned}\end{equation}
and
\begin{equation}
\sum_{k=1}^\infty\f{\bi{2k}k16^k(46k^2-11k+1)(292H_{2k-1}^{(2)}-77H_{k-1}^{(2)})}{k^2(2k-1)^2\bi{3k}k\bi{6k}{3k}}
=\f{178}3\pi^4.
\end{equation}
\end{conjecture}
\begin{remark} Note that \eqref{H631}$-6\times$\eqref{H2k-1} yields the identity
\begin{align*}&\sum_{k=1}^\infty\f{\bi{2k}k16^k((46k^2-11k+1)(2H_{6k-1}-H_{3k-1}-6H_{2k-1}+H_{k-1})-14k)}{k^2(2k-1)^2\bi{3k}k\bi{6k}{3k}}
\\&\qquad=-56\zeta(3).\end{align*}
\end{remark}

\begin{conjecture}[2023-10-18] {\rm (i)} We have
\begin{equation}\label{64^kC(2k,k)}\sum_{k=1}^\infty\f{(22k^2-7k+1)64^k\bi{2k}k}{k^2(2k-1)^2\bi{3k}k\bi{6k}{3k}}=4\pi^2.
\end{equation}
Also,
\begin{equation}\begin{aligned}
&\sum_{k=1}^\infty\f{\bi{2k}k64^k((22k^2-7k+1)H_{k-1}+P(k)/(5k(2k-1)))}{k^2(2k-1)^2\bi{3k}k\bi{6k}{3k}}
\\&\qquad=\f{24}5(\pi^2\log2+14\zeta(3))\end{aligned}
\end{equation}
and
\begin{equation}\begin{aligned}
&\sum_{k=1}^\infty\f{\bi{2k}k64^k((22k^2-7k+1)H_{2k-1}+Q(k)/(5k(2k-1)))}{k^2(2k-1)^2\bi{3k}k\bi{6k}{3k}}
\\&\qquad=\f{16}5(3\pi^2\log2+7\zeta(3)),\end{aligned}
\end{equation}
where
$$P(k)=296k^3-60k^2+6k-1\ \t{and}\ Q(k)=142k^3-45k^2+12k-2.$$
Moreover,
\begin{equation}\sum_{k=1}^\infty\f{\bi{2k}k64^k((22k^2-7k+1)\mathcal H(k)-8k(3k-1)/(2k-1))}
{k^2(2k-1)^2\bi{3k}k\bi{6k}{3k}}=0,
\end{equation}
where $\mathcal H(k)=2H_{6k-1}-H_{3k-1}-2H_{2k-1}$.

{\rm (ii)} We have
\begin{equation}\sum_{k=1}^\infty\f{\bi{2k}k64^k\l((22k^2-7k+1)(H_{2k-1}^{(2)}-\f3{16}H_{k-1}^{(2)})
-\f{3k(6k-1)}{(2k-1)^2}\r)}{k^2(2k-1)^2\bi{3k}k\bi{6k}{3k}}
=\f{\pi^4}{12}.
\end{equation}
\end{conjecture}
\begin{remark} The series in \eqref{64^kC(2k,k)} has converging rate $16/27$.
\end{remark}

\begin{conjecture} [2023-08-23]
 We have
\begin{equation}\label{16^kP(k)}
\sum_{k=1}^\infty\f{\bi{2k}k16^kf(k)}{k^2(2k-1)^2(6k-1)(6k-5)\bi{3k}k\bi{6k}{3k}}=\pi^2,
\end{equation}
where $f(k):=276k^3-248k^2+69k-5.$
Also,
\begin{equation}\begin{aligned}&\sum_{k=1}^\infty\f{\bi{2k}k16^k(f(k)(8H_{6k-1}-4H_{3k-1}+H_{2k-1}-H_{k-1})-f_1(k))}
{k^2(2k-1)^2(6k-1)(6k-5)\bi{3k}k\bi{6k}{3k}}
\\&\qquad=\f 92\pi^2+10\pi^2\log2
\end{aligned}
\end{equation}
and
\begin{equation}\begin{aligned}&\sum_{k=1}^\infty\f{\bi{2k}k16^k(f(k)(2H_{6k-1}-H_{3k-1}-H_{2k-1}+H_{k-1})-f_2(k))}
{k^2(2k-1)^2(6k-1)(6k-5)\bi{3k}k\bi{6k}{3k}}
=3\pi^2,
\end{aligned}
\end{equation}
where
$$f_1(k):=54k^2+222k-\f{235}2+\f{25}{2k}\ \t{and}
\ f_2(k)=186k^2-227k-15-\f{30}{2k-1}.$$
\end{conjecture}
\begin{remark} The series in \eqref{16^kP(k)} has converging rate $4/27$.
\end{remark}

\begin{conjecture}[2023-08-23] Let $P(k)=828k^3-888k^2+207k-11$.

{\rm (i)} We have
\begin{equation}\label{P(k)}\sum_{k=1}^\infty\f{\bi{2k}k16^kP(k)}{k(2k-1)^2(6k-1)(6k-5)\bi{3k}k\bi{6k}{3k}}=\f32\pi^2,
\end{equation}
\begin{equation}\begin{aligned}
&\sum_{k=1}^\infty\f{\bi{2k}k16^k(P(k)(H_{2k-1}-H_{k-1})-(6k-1)^2(138k-109)/(4(2k-1)))}{k(2k-1)^2(6k-1)(6k-5)\bi{3k}k\bi{6k}{3k}}
\\&\qquad=3\pi^2\log2-21\zeta(3),
\end{aligned}\end{equation}
and
\begin{equation}\begin{aligned}
&\sum_{k=1}^\infty\f{\bi{2k}k16^k(P(k)(2H_{6k-1}-H_{3k-1})+(6k-1)(84k^2-4k-79)/(4(2k-1)))}{k(2k-1)^2(6k-1)(6k-5)\bi{3k}k\bi{6k}{3k}}
\\&\qquad=3\pi^2\log2+21\zeta(3).
\end{aligned}\end{equation}

{\rm (ii)} We have
\begin{equation}\begin{aligned}
&\sum_{k=1}^\infty\f{\bi{2k}k16^k(P(k)(5H_{2k-1}^{(2)}-H_{k-1}^{(2)})-(6k-1)^2(231k-185)/(2k-1)^2)}{k(2k-1)^2(6k-1)(6k-5)\bi{3k}k\bi{6k}{3k}}
\\&\qquad=-\f{\pi^4}2
\end{aligned}\end{equation}
and
\begin{equation}\begin{aligned}
&\sum_{k=1}^\infty\f{\bi{2k}k16^k(P(k)(556H_{6k-1}^{(2)}-139H_{3k-1}^{(2)}-692H_{2k-1}^{(2)}+194H_{k-1}^{(2)})+4g(k))}{k(2k-1)^2(6k-1)(6k-5)\bi{3k}k\bi{6k}{3k}}
\\&\qquad=\f{55}2\pi^4,
\end{aligned}\end{equation}
where $$g(k):=\f{6090k^2+4225k-2217}{2k-1}.$$
\end{conjecture}
\begin{remark} The series in \eqref{P(k)} has converging rate $4/27$.
\end{remark}
\begin{conjecture} [2023-08-23] We have
\begin{equation}\sum_{k=1}^\infty\f{\bi{2k}k16^k(828k^3+1320k^2-745k+65)}{k^2(2k-1)(6k-1)(6k-5)\bi{3k}k\bi{6k}{3k}}=16\pi^2.
\end{equation}
\end{conjecture}
\begin{remark} We haven't found any variant of this identity with summands involving harmonic numbers.
\end{remark}


\begin{thebibliography}{99}

\bibitem{AKP}
\newblock G. Almkvist, C. Krattenthaler and J. Petersson,
\newblock {\it Some new series for $\pi$},
\newblock {Experiment. Math.}  \textbf{12} (2003), 441--456.

\bibitem{Ba} N. Batir, {\it On the series $\sum_{k=1}^\infty\bi{3k}k^{-1}k^{-n}x^k$},
Proc. Indian Acad. Sci. (Math. Sci.) {\bf 115} (2005), 371--381.

\bibitem{Be} F. Bellard, $\pi$ Formulas, Algorithms and Computations, {\tt https://bellard.org/pi}.

\bibitem{BC} J. M. Borwein and M. Chamberland,
 {\it Integer powers of arcsin},
\newblock {Int. J. Math. Math. Sci.}
\textbf{2007} (2007), Article ID 19381, 10pp.

\bibitem{CZ} W. Chu and W. Zhang, {\it Accelerating Dougall's ${}_5F_4$-sum and infinite series involving $\pi$}, Math. Comp. {\bf 83} (2014), 475--512.


\bibitem{Co} S. Cooper, {Ramanujan's Theta Functions}, Springer, Cham, 2017.

\bibitem{Wei} C. Wei, {\it On two double series for $\pi$ and their $q$-analogues},
Ramanujan J. {\bf 60} (2023), 615--625.


\bibitem{S13} Z.-W. Sun, {\it Products and sums divisible by central binomial coefficients},
 Electron. J. Combin.  {\bf 20} (2013), no.\,1, \#P9, 1-14.

 \bibitem{BNK} Z.-W. Sun, {\it New congruences involving harmonic numbers}, Nanjing Univ. J. Math.
 Biquarterly {\bf 40} (2023), 1--33.

 \bibitem{harmonic} Z.-W. Sun, {\it Series with summands involving harmonic numbers}, in: M. B. Nathanson (ed.),
 Combinatorial and Additive Number Theory VI, Springer, to appear. See also {\tt arXiv:2210.07238}.


\end{thebibliography}
\end{document}